\documentclass{amsart}
\usepackage{amsmath,amssymb,amsthm, mathtools,url}
\usepackage{physics}
\usepackage[margin=1.16in]{geometry}
\usepackage{tikz}
\usetikzlibrary{shapes.geometric, arrows}
\tikzstyle{arrow} = [thick,->,>=stealth]
\usetikzlibrary{hobby}
\usepackage{amsmath}
\usepackage[font = small, width=.6\linewidth]{caption}
\usepackage{amsthm}
\usepackage{amsfonts}
\usepackage{enumerate}
\definecolor{carnelian}{rgb}{0.7, 0.11, 0.11}
\definecolor{greenish}{rgb}{0.7, 1, 0.6}
\definecolor{bubblegum}{rgb}{0.99, 0.76, 0.8}
\definecolor{bittersweet}{rgb}{1.0, 0.44, 0.37}
\definecolor{babypink}{rgb}{0.96, 0.76, 0.76}
\definecolor{carnationpink}{rgb}{1.0, 0.65, 0.79}
\definecolor{cosmiclatte}{rgb}{1.0, 0.97, 0.91}
\definecolor{bisque}{rgb}{1.0, 0.89, 0.77}
\definecolor{celadon}{rgb}{0.67, 0.88, 0.69}
\usepackage{upgreek}
\usetikzlibrary{shapes.geometric, arrows.meta, positioning, decorations.pathreplacing}
\numberwithin{equation}{section}

\newtheorem{theorem}{Theorem}[section]
\newtheorem{corollary}[theorem]{Corollary}
\newtheorem{proposition}[theorem]{Proposition}
\newtheorem{definition}[theorem]{Definition}

\newtheorem{lemma}[theorem]{Lemma}

\newtheorem*{thmA}{Theorem A}
\newtheorem*{thmB}{Theorem B}

\newtheorem{strat}[theorem]{Strategy}
\newtheorem{ppp}[theorem]{Problem}

\theoremstyle{definition}
\newtheorem{example}[equation]{Example}

\title[Hidden ideals in pointed posets]{{\bf An efficient search strategy for hidden ideals\\ in pointed partially ordered sets}}

\author{\sc Roma Eisel}
\address{Duquesne University\\ Pittsburgh\\ PA, USA 15282}
\email{eiselr@duq.edu}

\author{\sc Valerie McMullen}
\address{Duquesne University\\ Pittsburgh\\ PA, USA 15282}
\email{mcmullenv@duq.edu}

\author{\sc Robert Muth}
\address{Department of Mathematics and Computer Science\\ Duquesne University\\ Pittsburgh\\ PA, USA 15282}
\email{muthr@duq.edu}

\begin{document}

\begin{abstract}
We consider a combinatorial question about searching for an unknown ideal \(\mu\) within a known pointed poset \(\lambda\). Elements of \(\lambda\) may be queried for membership in \(\mu\), but at most \(k\) positive queries are permitted. We provide a general search strategy for this problem, and establish new bounds (based on \(k\) and the  degree and height of \(\lambda\)) for the total number of queries required to identify \(\mu\). We show that this strategy performs asymptotically optimally on the family of complete \(\ell\)-ary trees as the height grows.
\end{abstract}

\maketitle

\section{Introduction}
\subsection{The hidden ideal problem}\label{HIprob}
Let \(\lambda\) be a finite partially ordered set, and let \(k \in \mathbb{N}\). Hidden within \(\lambda\) is a (possibly empty) unknown ideal \(\mu\). We seek to identify  \(\mu\) by sequentially querying elements of \(\lambda\) for membership in \(\mu\), under the restriction that at most \(k\) positive query results are permitted. Letting \(\boldsymbol{q}_k(\lambda)\) represent the minimum total number of queries needed to guarantee identification of \(\mu\), the {\em hidden ideal problem} is to solve:

\begin{ppp}\label{P1}
Find \(\boldsymbol{q}_k(\lambda)\), and identify a search strategy which realizes this value.
\end{ppp}

This problem has been studied under numerous guises and for various families of posets.
When \(\lambda\) is totally ordered, Problem~\ref{P1} is related to the `\(k\)-egg' or `\(k\)-marble' problem \cite{marbles, Egg1, Egg3}, which is featured in numerous texts on dynamic programming and optimization. When \(\lambda\) is a rectangular poset with one dimension not more than six, Problem~\ref{P1} was settled in \cite{IJM}. 
More generally, this problem relates to a broad body of work in computer science and optimization on efficient search within sets with partial order---see for instance \cite{Ben, Chen, Cic2, Cic3, Cic, Sort, Deren}---motivated by such applications as debugging, file synchronization, and information retrieval. From this point of view Problem~\ref{P1} is a consideration of this topic under an additional restriction, or a `cost' imposition on excessive positive queries. 

\subsection{A general search strategy}\label{gensearsec}
Contrary to \cite{marbles, IJM, Egg1, Egg3}, which focus on {\em optimal} search strategies for {\em specific} families of posets, our goal in the present paper is to produce {\em efficient} search strategies for {\em general} pointed posets (i.e., those with a unique minimal element). Generally speaking, the only known upper bound for \(\boldsymbol{q}_k(\lambda)\) is given by \(|\lambda|\)---a bound achieved via the naive strategy of querying all nodes of \(\lambda\) in a top-down fashion. 

In Strategy~\ref{thestrat} we present a new query strategy which may be applied to any pointed poset \(\lambda\) (in Example~\ref{exsec} we provide a short example of Strategy~\ref{thestrat} at work). We then prove the following main result, which provides a bound on the number of total queries required based on the {\em degree} and {\em height} of the poset \(\lambda\) (see \S\ref{prelimsec}).

\begin{thmA}
Let \(k \in \mathbb{N}\), and let \(\lambda\) be a pointed poset of degree \(\ell\) and height \(n\). Then Strategy~\ref{thestrat} guarantees identification of the hidden ideal \(\mu\), using at most \(k\) positive queries, and at most
\begin{align*}
f_{k,\ell}(n) := \sum_{i=0}^{\lceil n/k \rceil-1} \ell^i + \sum_{j=1}^{k-1} \ell^{\lceil (n-j)/k\rceil}
\leq
\begin{cases}
\left\lceil \frac{n+1}{k}\right \rceil  &\textup{if } \ell=1;\\
k\ell^{\lceil n/k \rceil} &\textup{if }\ell > 1
\end{cases}
\end{align*}
total queries.
\end{thmA}

This appears as Theorem~\ref{mainthmA} in the body. It follows from Theorem A then that \(\boldsymbol{q}_k(\lambda) \leq f_{k,\ell}(n)\) for any pointed poset \(\lambda\) of degree \(\ell\) and height \(n\), yielding a bound that greatly improves on the old bound \(\boldsymbol{q}_k(\lambda) \leq |\lambda|\) in general. For instance, if \(\lambda\) is the complete \(\ell\)-ary tree of height \(n\) (see \S\ref{treesec}), then \(\boldsymbol{q}_k(\lambda) \leq f_{k,\ell}(n) \leq k \ell^{(k+1)/k}|\lambda|^{1/k}\).

\subsection{Asymptotic optimality}
It was shown in \cite[Proposition 2.2]{IJM} that \(\boldsymbol{q}_k(\lambda)\) has a recursive combinatorial description:
\begin{align*}
\boldsymbol{q}_k(\lambda) = 
\begin{cases}
0 &\textup{if }\lambda = \varnothing;\\
|\lambda| &\textup{if }k = 1;\\
\min \{ \max\{ \boldsymbol{q}_k(\lambda_{\not \succcurlyeq u}) , \boldsymbol{q}_{k-1}(\lambda_{\succ u})\}\mid u \in \lambda\}+1
&  \textup{if }k >1, \lambda \neq \varnothing.
\end{cases}
\end{align*}
For reasons of storage and computation time, it is quite demanding to compute \(\boldsymbol{q}_k(\lambda)\) via this formula for large posets, and optimal search strategies as in Problem~\ref{P1} are linked to the specific structure of \(\lambda\) in very delicate ways (see for instance \cite[Algorithm 6.2]{IJM}). As Strategy~\ref{thestrat} is a `one-size-fits-all' approach on the family of pointed posets of a given degree and height, it generally does not produce solutions within an optimally minimal number of total queries, but does in a specific sense perform well on large members of this family, as we now explain.

For \(\ell, n \in \mathbb{N}\), let \(\mathcal{T}_{\ell}(n)\) be the complete \(\ell\)-ary tree of height \(n\). Then \(\mathcal{T}_{\ell}(n)\) has maximal cardinality among the family of pointed posets of degree \(\ell\) and height \(n\). Our second main result shows that Strategy~\ref{thestrat} performs {\em asymptotically optimally} on complete \(\ell\)-ary trees as the height \(n\) increases.

\begin{thmB}
Fix \(k,\ell \in \mathbb{N}\).
Then there exist \(m,M > 0\) such that
\begin{align*}
m \ell^{\,n/k} \leq \boldsymbol{q}_{k}(\mathcal{T}_\ell(n)) \leq f_{k,\ell}(n) \leq M \ell^{\,n/k},
\end{align*}
for all \(n \in \mathbb{N}\). 
Thus \(\boldsymbol{q}_{k}(\mathcal{T}_\ell(n)) = \Theta(\ell^{\,n/k}) = f_{k,\ell}(n)\) as functions of \(n\).
\end{thmB}

This appears as Theorem~\ref{growthm} in the body. 
More broadly, this speaks to the overall efficiency of Strategy~\ref{thestrat}. If we fix a poset degree \(\ell\) and a limit \(k\) on the number of positive queries, then Theorems A and B guarantee that Strategy~\ref{thestrat} will identify a hidden ideal \(\mu\) in a pointed poset \(\lambda\) of height \(n\) in \(\mathcal{O}(\ell^{n/k})\) total queries. On the other hand, Theorem B shows there always exist such posets (e.g. the complete \(\ell\)-ary trees) which {\em require} \(\Omega(\ell^{n/k})\) total queries to solve.

\subsection{Acknowledgements}
This paper was initially motivated by a problem introduced by Jonathan Kujawa in \cite{Kuj}, which recast the search for pollution sources in the Mississippi river system as a hidden ideal problem in a tree poset. We thank him for the inspiration.

\section{Preliminaries}\label{prelimsec}
In this section we briefly lay out preliminaries and relevant definitions. Throughout, we take \(\mathbb{N} = \mathbb{Z}_{>0}\), and use the notation \([a,b] := \{a,a+1,\ldots, b\}\) for any integers \(a\leq b\). 

\subsection{Posets and ideals}
We give now a very brief introduction to posets and related terms. 
The reader should see \cite{Davey, Dushnik} for a more thorough treatment of the subject.
    \begin{definition}
        A \textit{partially ordered set}, or \textit{poset}, is a set $\lambda$ together with a binary relation \(\succcurlyeq\), which satisfies the following conditions for all $u, v, w \in \lambda$:
    \begin{enumerate}[(i)]
                \item $u \succcurlyeq u$ (reflexivity);
                \item $u \succcurlyeq v$ and $v \succcurlyeq u$ implies $u = v$ (antisymmetry);
                \item $u \succcurlyeq v$ and $v \succcurlyeq w$ implies $u \succcurlyeq w$ (transitivity).
    \end{enumerate} 
    \end{definition}
From this point forward, we assume all posets under consideration are {\em finite}. If \(u \succcurlyeq v\) then we say that \(u\) {\em dominates} \(v\).
    For \(u, v \in \lambda\), we write \(u \succ v\) provided \(u \succcurlyeq v\) and \(u \neq v\). We say a poset \(\lambda\) is {\em pointed} provided that for all \(a,b \in \lambda\), there exists \(c \in \lambda\) with \(a \succcurlyeq c\) and \(b \succcurlyeq c\). In our context of finite posets, {\em pointed} implies that either \(\lambda = \varnothing\) or \(\lambda\) has a unique minimal element. We will say two posets \(\lambda, \nu\) are {\em isomorphic}, and write \(\lambda \cong \nu\), if there exist order-preserving mutually inverse maps \( \lambda \rightleftarrows \nu\).    

\subsubsection{Subposets} 
If \(\mu\) is a subset of a poset \(\lambda\), then \(\mu\) is itself a poset under the partial order
inherited from \(\lambda\), and we always assume we take this partial order on \(\mu\). For \(v \in \lambda\), we define some particular subposets of \(\lambda\):
\begin{align*}
\lambda_{\succcurlyeq v} := \{u \in \lambda \mid u \succcurlyeq v\}
\qquad
\lambda_{\succ v} := \{u \in \lambda \mid u \succ v\} 
\qquad
\lambda_{\not\succcurlyeq v} := \{u \in \lambda \mid u \not\succcurlyeq v\},
\end{align*}
More generally, for \(S \subseteq \lambda\), we write:
\begin{align*}
\lambda_{\succcurlyeq S} &= \bigcup_{v \in S} \lambda_{\succcurlyeq v} = \{u \in \lambda \mid u \succcurlyeq v \textup{ for some }v \in S\}\\
\lambda_{\not\succcurlyeq S} &= \bigcap_{v \in S} \lambda_{\not\succcurlyeq v} = \{u \in \lambda \mid u \not\succcurlyeq v \textup{ for all }v \in S\} = \lambda - \lambda_{\succcurlyeq S}
\end{align*}
We similarly define \(\lambda_{\preccurlyeq v}, \lambda_{\preccurlyeq S}\), and so on. Note that \(\lambda_{\succcurlyeq v}\) is always a pointed poset with minimal element \(v\), and \(\lambda_{\not \succcurlyeq v}\) is pointed if and only if \(\lambda\) is pointed.

\subsubsection{Ideals}
\begin{definition}
Assume \(\mu \subseteq \lambda\). We say \(\mu\) is an {\em ideal} in \(\lambda\) provided that \(\mu = \varnothing\) or  \(\mu = \lambda_{\preccurlyeq v}\) for some \(v \in \lambda\). In the latter case we say \(\mu\) is {\em generated} by \(v\).
\end{definition}

In our context of finite sets, this notion of ideal is equivalent to the that of a {\em directed lower set}.

\subsubsection{Hasse diagrams}
        Let $\lambda$ be a poset, and let $u,v \in \lambda$. We say $u$ {\em covers} $v$ provided $u \succ v$ and there is no $w \in \lambda$ with $u \succ w \succ v$. We will visually display posets via their {\em Hasse diagram}---this is a graph whose vertices are the elements of \(\lambda\), with an edge drawn upward from \(v\) to \(u\) whenever \(u\) covers \(v\).

        \begin{example}\label{firstposetex} In Figure~\ref{hidden_ideal_orange} we display the Hasse diagram of a pointed poset \(\lambda\). The unique minimal element \(m \in \lambda\) is shown. The ideal \(\lambda_{\preccurlyeq v}\) generated by \(v \in \lambda\) is also shown.
        \end{example}
        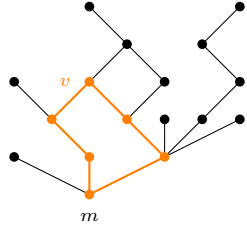
\begin{figure}[h]
        \begin{align*}
            \begin{tikzpicture}
                \draw (0,0)--(-1,.5);
                \draw [thick, orange](0,0)--(0,.5);
                \draw [thick, orange](0,0)--(1,.5);
                \draw [thick, orange](-.5,1)--(0,.5);
                \draw [thick, orange](1,.5)--(.5,1);
                \draw (1,.5)--(1,1);
                \draw (1,.5)--(1.5,1);
                \draw (1,.5)--(2,1);
                \draw (1.5,1)--(2,1.5);
                \draw (2,1.5)--(1.5,2);
                \draw (-.5,1)--(-1,1.5);
                \draw [thick, orange](-.5,1)--(0,1.5);
                \draw (0,1.5)--(0.5,2);
                \draw (0.5,2)--(0,2.5);
                \draw (0.5,2)--(1,1.5);
                \draw [thick, orange](0, 1.5)--(.5,1);
                \draw (.5,1)--(1,1.5);
                \draw (1.5,2)--(2,2.5);
                \draw[thick, orange, fill = orange](0,0) circle(.5mm);
                \draw[thick, black, fill = black](-1,.5) circle(.5mm);
                \draw[thick, orange, fill = orange](0,.5) circle(.5mm);
                \draw[thick, orange, fill = orange](1,.5) circle(.5mm);
                \draw[thick, orange, fill = orange](-.5,1) circle(.5mm);
                \draw[thick, orange, fill = orange](.5,1) circle(.5mm);
                \draw[thick, black, fill = black](1,1) circle(.5mm);
                \draw[thick, black, fill = black](1.5,1) circle(.5mm);
                \draw[thick, black, fill = black](2,1) circle(.5mm);
                \draw[thick, black, fill = black](2,1.5) circle(.5mm);
                \draw[thick, black, fill = black](1.5,2) circle(.5mm);
                \draw[thick, black, fill = black](2,2.5) circle(.5mm);
                \draw[thick, black, fill = black](-1,1.5) circle(.5mm);
                \draw[thick, orange, fill = orange](0,1.5) circle(.5mm);
                \draw[thick, black, fill = black](0.5,2) circle(.5mm);
                \draw[thick, black, fill = black](0,2.5) circle(.5mm);
                \draw[thick, black, fill = black](1,1.5) circle(.5mm);
                \node [label=left:\textcolor{orange}{\scriptsize $v$}] at (0,1.5) {};
                \node [label=below:\textcolor{black}{\scriptsize $m$}] at (0,0) {};
            \end{tikzpicture}
            \end{align*}
            \caption{ A pointed poset $\lambda$ with 17 elements and ideal $ \lambda_{\preccurlyeq v}$ highlighted in orange.}
                        \label{hidden_ideal_orange}
            \end{figure}

    \subsection{Height and degree}
    Now we introduce some important statistics on posets.

    \begin{definition}
    Let \(\lambda\) be a poset. 
    \begin{enumerate}[(i)]
    \item For \(v \in \lambda\), we define the {\em degree} of \(v\) to be
    \begin{align*}
    \deg(v) := \# \{u \in \lambda \mid u \textup{ covers } v\}.
    \end{align*}
    \item We define the {\em degree} of \(\lambda\) to be
    \begin{align*}
    \deg(\lambda) := \max\{\deg(v) \mid v \in \lambda\}.
    \end{align*}
    \end{enumerate}
    \end{definition}

\begin{example}
For the poset \(\lambda\) in Figure~\ref{hidden_ideal_orange}, we have \(\deg(v) = 1\), \(\deg(m) = 3\), and \(\deg(\lambda) = 4\).
\end{example}

A subset $C$ of a pointed poset \(\lambda\) is a {\em chain} provided every pair of nodes $x, y \in C$ is comparable. We will label the nodes of \(C\) as 
\begin{align*}
C = \{C_{|C|} \succ C_{|C| - 1} \succ \cdots \succ C_1\}
\end{align*}
where \(C_1, \ldots, C_{|C|} \in \lambda\). 
For \(v \in \lambda\), we say \(C\) is a {\em \(v\)-chain} if \(C_{|C|} = v\). We say moreover \(C\) is a {\em maximal} \(v\)-chain if \(|C| \geq |D|\) for every \(v\)-chain \(D\). For a maximal \(v\)-chain \(C\), we have that \(C_{t}\) covers \(C_{t-1}\) for \(t \in [2,|C|]\), and \(C_1\) is the unique minimal element of the pointed poset \(\lambda\).

Finally, we say a chain \(C\) is a \textit{maximal} \(\lambda\)-chain provided that \(|C| \geq |D|\) for every chain \(D\) in \(\lambda\).

    \begin{definition}
        Let $\lambda$ be a poset. 
        \begin{enumerate}[(i)]
        \item For \(v \in \lambda\), we define the {\em height} of \(v\) to be \(\textup{ht}(v) := |C|\), where \(C\) is any maximal \(v\)-chain.

        \item We define the {\em height} of \(\lambda\) to be \(\textup{ht}(\lambda) := |C|\), where \(C\) is any maximal \(\lambda\)-chain.

        \end{enumerate}
        \end{definition}

        We make a few remarks on height, \(v\)-chains and \(\lambda\)-chains that are quickly verified. 
        First, note that \(\textup{ht}(\lambda) = \max\{\textup{ht}(v) \mid v \in \lambda\}\), and \(\textup{ht}(\lambda_{\preccurlyeq v}) = \textup{ht}(v)\) for every \(v \in \lambda\). If \(C\) is a maximal \(\lambda\)-chain, or maximal \(v\)-chain for some \(v \in \lambda\), then \(\textup{ht}(C_t) = t\) for every \(t \in [1,|C|]\).

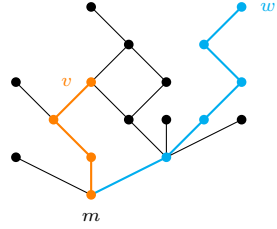
\begin{figure}[h]
        \begin{align*}
            \begin{tikzpicture}
                \draw (0,0)--(-1,.5);
                \draw [thick, orange](0,0)--(0,.5);
                \draw [thick, cyan](0,0)--(1,.5);
                \draw [thick, orange](-.5,1)--(0,.5);
                \draw (1,.5)--(.5,1);
                \draw (1,.5)--(1,1);
                \draw [thick, cyan] (1,.5)--(1.5,1);
                \draw (1,.5)--(2,1);
                \draw [thick, cyan] (1.5,1)--(2,1.5);
                \draw [thick, cyan] (2,1.5)--(1.5,2);
                \draw (-.5,1)--(-1,1.5);
                \draw [thick, orange](-.5,1)--(0,1.5);
                \draw (0,1.5)--(0.5,2);
                \draw (0.5,2)--(0,2.5);
                \draw (0.5,2)--(1,1.5);
                \draw (0, 1.5)--(.5,1);
                \draw (.5,1)--(1,1.5);
                \draw [thick, cyan](1.5,2)--(2,2.5);
                \draw[thick, orange, fill = orange](0,0) circle(.5mm);
                \draw[thick, black, fill = black](-1,.5) circle(.5mm);
                \draw[thick, orange, fill = orange](0,.5) circle(.5mm);
                \draw[thick, cyan, fill = cyan](1,.5) circle(.5mm);
                \draw[thick, orange, fill = orange](-.5,1) circle(.5mm);
                \draw[thick, black, fill = black](.5,1) circle(.5mm);
                \draw[thick, black, fill = black](1,1) circle(.5mm);
                \draw[thick, cyan, fill = cyan](1.5,1) circle(.5mm);
                \draw[thick, black, fill = black](2,1) circle(.5mm);
                \draw[thick, cyan, fill = cyan](2,1.5) circle(.5mm);
                \draw[thick, cyan, fill = cyan](1.5,2) circle(.5mm);
                \draw[thick, cyan, fill = cyan](2,2.5) circle(.5mm);
                \draw[thick, black, fill = black](-1,1.5) circle(.5mm);
                \draw[thick, orange, fill = orange](0,1.5) circle(.5mm);
                \draw[thick, black, fill = black](0.5,2) circle(.5mm);
                \draw[thick, black, fill = black](0,2.5) circle(.5mm);
                \draw[thick, black, fill = black](1,1.5) circle(.5mm);
                \node [label=left:\textcolor{orange}{\scriptsize $v$}] at (0,1.5) {};
                \node [label=right:\textcolor{cyan}{\scriptsize $w$}] at (2,2.5) {};
                \node [label=below:\textcolor{black}{\scriptsize $m$}] at (0,0) {};
            \end{tikzpicture}
            \end{align*}
            \caption{ The pointed poset $\lambda$ with highlighted chains.}
                        \label{chainspic}
            \end{figure}

    \begin{example}
For the poset \(\lambda\) in Figure~\ref{chainspic}, we have highlighted two chains. The orange chain is a maximal \(v\)-chain. The cyan chain is a maximal \(w\)-chain, and a maximal \(\lambda\)-chain. We have
\begin{align*}
\textup{ht}(v) = 4 \qquad \textup{ht}(w) = 6 \qquad \textup{ht}(m) = 1 \qquad \textup{ht}(\lambda) = 6.
\end{align*}
\end{example}

    \subsection{The posets \(\mathcal{P}(\ell,n)\)}
    Going forward, we write 
    \begin{align*}
    \mathcal{P}(\ell, n) = \{ \lambda \textup{ a pointed poset} \mid \deg(\lambda) \leq \ell, \; \textup{ht}(\lambda) \leq n).
    \end{align*}

    For example, the poset \(\lambda\) of Figures~\ref{hidden_ideal_orange} and \ref{chainspic} is a member of \(\mathcal{P}(4,6)\). 

    \begin{lemma}\label{new_lem}
        Let $\lambda \in \mathcal{P}(\ell,n)$. Then \(\lambda\) has at most \(\ell^{t-1}\) nodes of height \(t \in [1,n]\), and $| \lambda | \leq \sum_{t=0}^{n-1} \ell^t$.
    \end{lemma}

    \begin{proof}
    For \(t \in [1,n]\), write \(m_t\) for the number of nodes in \(\lambda\) of height \(t\). Note that since \(\lambda\) is pointed, we must have \(m_1 = 1\). Since the degree of \(\lambda\) is \(\ell\), we also have that \(m_t \leq \ell m_{t-1}\) for \(t \in [2,n]\). Therefore it follows that \(m_t \leq \ell^{t-1}\) for \(t \in [1,n]\). Thus
    \begin{align}\label{mtsum}
    |\lambda| = \sum_{t=1}^n m_t \leq \sum_{t=1}^n \ell^{t-1} =  \sum_{t=0}^{n-1} \ell^{t},
    \end{align}
    as required.
    \end{proof}

    \subsection{Complete \(\ell\)-ary trees}\label{treesec}
    
    \begin{definition}
    We say a pointed poset \(\lambda \in \mathcal{P}(\ell,n)\) is a {\em complete \(\ell\)-ary tree of height \(n\)} provided that \(\lambda\) has exactly \(\ell^{t-1}\) nodes of height \(t \in [1,n]\). 
    \end{definition}

    Thus the complete \(\ell\)-ary trees of height \(n\) are exactly the posets in \(\mathcal{P}(\ell,n)\) which render the inequality in Lemma~\ref{new_lem} an equality. If \(\lambda \in \mathcal{P}(\ell,n)\) is a complete \(\ell\)-ary tree of height \(n\), it follows that every node of height less than \(n\) has degree \(\ell\), and \(\lambda_{\preccurlyeq v}\) is a chain for every \(v \in \lambda\). All complete \(\ell\)-ary trees of height \(n\) are isomorphic, and we use the notation \(\mathcal{T}_\ell(n)\) to indicate (the isomorphism class of) the complete \(\ell\)-ary tree of height \(n\).

\begin{figure}[h]
        \begin{align*}
            \begin{tikzpicture}
                \draw[thick, black, fill = black](8, 0) circle(.5mm);
                \draw[thick, black, fill = black](8, 1) circle(.5mm);
                \draw[thick, black, fill = black](5, 1) circle(.5mm);
                \draw[thick, black, fill = black](11, 1) circle(.5mm);
                \draw[thick, black, fill = black](8, 2) circle(.5mm);
                \draw[thick, black, fill = black](5, 2) circle(.5mm);
                \draw[thick, black, fill = black](11, 2) circle(.5mm);
                \draw[thick, black, fill = black](4, 2) circle(.5mm);
                \draw[thick, black, fill = black](6, 2) circle(.5mm);
                \draw[thick, black, fill = black](7, 2) circle(.5mm);
                \draw[thick, black, fill = black](9, 2) circle(.5mm);
                \draw[thick, black, fill = black](10, 2) circle(.5mm);
                \draw[thick, black, fill = black](12, 2) circle(.5mm);
                \draw[thick, black, fill = black](2, 1) circle(.5mm);
                \draw[thick, black, fill = black](1, 2) circle(.5mm);
                \draw[thick, black, fill = black](2, 2) circle(.5mm);
                \draw[thick, black, fill = black](3, 2) circle(.5mm);
                \draw[thick, black, fill = black](14, 1) circle(.5mm);
                \draw[thick, black, fill = black](13, 2) circle(.5mm);
                \draw[thick, black, fill = black](14, 2) circle(.5mm);
                \draw[thick, black, fill = black](15, 2) circle(.5mm);
                \draw (8, 0)--(2, 1);
                \draw (8, 0)--(8, 1);
                \draw (8, 0)--(5, 1);
                \draw (8, 0)--(11, 1);
                \draw (8, 0)--(14, 1);
                \draw (8, 1)--(7, 2);
                \draw (8, 1)--(8, 2);
                \draw (8, 1)--(9, 2);
                \draw (5, 1)--(4, 2);
                \draw (5, 1)--(6, 2);
                \draw (5, 1)--(5, 2);
                \draw (11, 1)--(10, 2);
                \draw (11, 1)--(11, 2);
                \draw (11, 1)--(12, 2);
                \draw (14, 1)--(13, 2);
                \draw (14, 1)--(14, 2);
                \draw (14, 1)--(15, 2);
                \draw (2, 1)--(1, 2);
                \draw (2, 1)--(2, 2);
                \draw (2, 1)--(3, 2);
\draw (2,1)--(1.5,2);
\draw (2,1)--(2.5,2);
\draw (5,1)--(4.5,2);
\draw (5,1)--(5.5,2);
\draw (8,1)--(7.5,2);
\draw (8,1)--(8.5,2);
\draw (11,1)--(10.5,2);
\draw (11,1)--(11.5,2);
\draw (14,1)--(13.5,2);
\draw (14,1)--(14.5,2);
\draw[thick, black, fill = black](1.5, 2) circle(.5mm);
\draw[thick, black, fill = black](2.5, 2) circle(.5mm);
\draw[thick, black, fill = black](4.5, 2) circle(.5mm);
\draw[thick, black, fill = black](5.5, 2) circle(.5mm);
\draw[thick, black, fill = black](7.5, 2) circle(.5mm);
\draw[thick, black, fill = black](8.5, 2) circle(.5mm);
\draw[thick, black, fill = black](10.5, 2) circle(.5mm);
\draw[thick, black, fill = black](11.5, 2) circle(.5mm);
\draw[thick, black, fill = black](13.5, 2) circle(.5mm);
\draw[thick, black, fill = black](14.5, 2) circle(.5mm);
            \end{tikzpicture}
            \end{align*}
\caption{\(\mathcal{T}_5(3)\), the complete \(5\)-ary tree of height 3.}
            \label{T53ex}
            \end{figure}
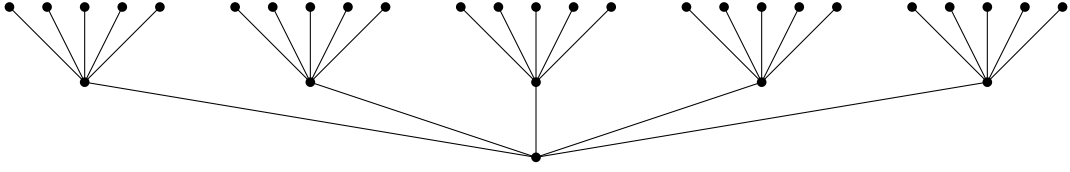

\begin{example}
In Figure~\ref{T53ex} we show \(\mathcal{T}_5(3)\), the complete \(5\)-ary tree of height \(3\).
\end{example}

\subsection{A useful technical lemma on posets}

    \begin{lemma} \label{lem_height<j}
        Let \(\lambda\) be a pointed poset, and assume $\lambda$ has exactly $t$ nodes $v_1, v_2, \dots, v_t$ of height $j$. Assume moreover that $w_1,w_2,\dots, w_t \in \lambda $ are such that: 
            \begin{itemize}
                \item $\textup{ht}(w_i) \leq j$, and;
                \item $w_i$ belongs to a maximal $\lambda_{\not\succcurlyeq\{w_1,w_2,\dots,w_{i-1}\}}$-chain,
            \end{itemize}
        for all $i \in [1,t]$. Then $\textup{ht}(\lambda_{\not\succcurlyeq\{w_1,\dots,w_t\}}) < j$.
    \end{lemma}

    \begin{proof}
        We go by induction on $t$. \\
        
        \noindent\textit{Base Case:} Assume $t = 1$. Then \(v_1\) is the only node of height \(j\) in \(\lambda\).  The node \(w_1\) is in some maximal \(\lambda\)-chain \(C\) by assumption, which contains one node of every height \(\leq \textup{ht}(\lambda)\), and thus contains a node of height \(j\). Thus \(v_1 \in C\). Since \(\textup{ht}(w_1) \leq j\), this implies that \(w_1 \preccurlyeq v_1\). 
        
        Now assume that \(u \in \lambda\) has \(\textup{ht}(u)\geq j\). Let \(D\) be a maximal \(u\)-chain in \(\lambda\). Then since \(D\)
        contains a node of height \(j\) as well, and hence contains \(v_1\). Since \(\textup{ht}(u)\geq \textup{ht}(v_1)\), this implies that \(u \succcurlyeq v_1 \succcurlyeq w_1 \). Therefore \(u \notin \lambda_{\not\succcurlyeq w_1}\). Thus \(\lambda_{\not\succcurlyeq w_1}\) contains only nodes of \(\lambda\) which have height less than \(j\), so \(\textup{ht}(\lambda_{\not \succcurlyeq w_1}) < j\), as required.\\

        \noindent\textit{Induction Step:} 
        Now let $t > 1$, and assume the claim holds for any $t' < t$. Again, by assumption, $w_1$ is in a maximal chain \(C\) in $\lambda$. This chain must include some height $j$ node $v_i$. Then, $v_i \succcurlyeq w_1$. Let $S \subseteq [1,t]$ be the set of all \(s\) such that $v_s \succcurlyeq w_1$ for all $s \in S$, and $S \neq \varnothing$. Then, the set of vertices of height $j$ in $\lambda_{\not\preccurlyeq\{w_1\}}$ is $\{v_k \mid k \not\in S\}$, and there are $t'=t - |S|$ of them. Set $\lambda' = \lambda_{\not\preccurlyeq\{w_1\}}$, and let $v'_1,\dots,v'_{t'}$ be the set of nodes of height $j$ in $\lambda'$. Then, $t' < t$. We also set $w'_b = w_{b+1}$ for $b \in [ 1,t']$. By induction, $\textup{ht}(\lambda'_{\not\succcurlyeq\{w'_1,\dots,w'_{t'}\}}) < j$. But we have
        \begin{align*}
        \lambda'_{\not\succcurlyeq\{w'_1,\dots,w'_{t'}\}} = (\lambda_{\not\succcurlyeq \{w_1\}})_{\not\succcurlyeq\{w_2,\dots,w_{t'+1}\}} \supseteq (\lambda_{\not\succcurlyeq\{w_1\}})_{\not\succcurlyeq\{w_2,\dots,w_t\}} = \lambda_{\not \succcurlyeq\{w_1, \ldots, w_t\}}.
        \end{align*}
        Thus, $\textup{ht}(\lambda_{\succcurlyeq\{w_1,\dots,w_t\}}) < j$, completing the induction step and the proof.
    \end{proof}

\section{A general strategy for the hidden ideal problem}
Recall now the {\em hidden ideal problem}, introduced in \S\ref{HIprob}. We describe in this section a general recursive strategy for identifying a hidden ideal in a pointed poset using a limited amount of positive queries. In \S\ref{stratanalsec} we establish bounds on the number of total queries this strategy entails, and in \S\ref{efficsec} we examine the asymptotic behavior of this strategy.

   \begin{strat}\label{thestrat}
Let \(\lambda^0\) be a pointed poset with a hidden ideal \(\mu\), and let \(k^0 \in \mathbb{N}\) be a limit on the number of positive queries permitted. 
Set \(\lambda:=\lambda^0\), \(k:=k^0\), and begin with the green block below. The strategy terminates by identifying \(\mu\) in one of the red blocks.
   \end{strat}
    \begin{align*}
    \begin{tikzpicture} [scale = .55]
        \node at (0,1) [rectangle, rounded corners, minimum height = 1cm, minimum width = 3cm, fill=celadon, draw] (1) {Input $\lambda$ and $k$};
%
        \node at (-7.5, 1) [rectangle, rounded corners, text centered, minimum height = 1cm, minimum width = 1.5cm, fill=bittersweet, draw] (2) {$\mu = \varnothing$};
%
        \node at (-6, -4.5) [rectangle, rounded corners, minimum height = 2cm, minimum width = 3cm, text width = 3cm, text centered, fill=bisque, draw] (3) {$\lambda = \{v\}$ for some $v$. Query $v$.};
            \node at (-8.5, -10) [rectangle, rounded corners, text centered, minimum height = 1cm, minimum width = 2cm, fill=bittersweet, draw] (6) {$\mu = \lambda^0_{\preccurlyeq v}$};
            \node at (-3.5, -10) [rectangle, rounded corners, text centered, minimum height = 1cm, minimum width = 2cm, fill=bittersweet, draw] (7) {$\mu = \varnothing$};
%
        \node at (6.5, -4.5) [rectangle, rounded corners, minimum height = 2cm, minimum width = 3cm, text width = 3cm, text centered, fill=bisque, draw] (4) {Set $n := \textup{ht}(\lambda)$};  
            \node at (14.2, -4.5) [rectangle, rounded corners, text centered, minimum height = 2cm, minimum width = 2cm, text width = 2cm, fill=bisque, draw] (5) {Set $v$ to be a node of height $n$};  
            \node at (3, -11) [rectangle, rounded corners, text centered, minimum height = 2cm, minimum width = 2cm, text width = 4cm, fill=bisque, draw] (8) {Set $v$ to be a node of height $2$ belonging to a maximal chain in $\lambda$}; 
            \node at (11, -11) [rectangle, rounded corners, text centered, minimum height =2cm, minimum width =4cm, text width = 3cm, fill=bisque, draw] (9) {Set $v$ to be a node of height $\left\lceil \frac{n+1}{k} \right\rceil$ belonging to a maximal chain in $\lambda$}; 
        \node at (6.5, -15.5) [rectangle, rounded corners, text centered, minimum height = 2cm, minimum width = 3cm, text width = 3cm, fill=bisque, draw] (10) {Query $v$.}; 
            \node at (1.5, -21.5) [rectangle, rounded corners, text centered, minimum height = 1.5cm, minimum width = 2cm, text width = 3cm, fill=bisque, draw] (11) {Set $\lambda := \lambda_{\succcurlyeq v}$ Set $k := k-1$};
            \node at (11.5, -21.5) [rectangle, rounded corners, text centered, minimum height = 1.5cm, minimum width = 2cm, text width = 3cm, fill=bisque, draw] (12) {Set $\lambda := \lambda_{\not\succcurlyeq v}$}; 
%
        \draw [arrow, thick] (1) -- node[rectangle, midway, black, fill = white] {\scriptsize $|\lambda | = 0$} (2);
        \draw [arrow, thick] (4) -- node[midway, fill = white] {\scriptsize $k = 1$} (5);
%
        \draw [arrow, thick] (12) -- (17.25,-21.5) -- (17.25,3) -- (1,3)--(1,1.9);
%
        \draw [thick] (-1.46, -22.15) -- node[midway, fill = white] {\scriptsize $|\lambda | > 1$} (-11,-22.15);
        \draw [thick] (-11,-22.15) -- (-11,3);
        \draw [thick] (-11,3) -- (-1,3);
        \draw [arrow, thick] (-1,3) -- (-1,1.9);
%
        \draw [arrow, thick] (1) -- node[midway, fill = white, sloped] {\scriptsize $|\lambda | = 1$} (3); 
        \draw [arrow, thick] (1) -- node[midway, fill = white, sloped] {\scriptsize $|\lambda | > 1$} (4);
%
        \draw [arrow, thick] (3) -- node[midway, fill = white, sloped] {\scriptsize $v \in \mu$} (6);
        \draw [arrow, thick] (3) -- node[midway, fill = white, sloped] {\scriptsize $v \notin \mu$} (7);
%
        \draw [arrow, thick] (4) -- node[midway, fill = white, sloped] {\scriptsize $k \geq n$} (8);
        \draw [arrow, thick] (4) -- node[midway, fill = white, sloped] {\scriptsize $1 \hspace{-0.35mm}<\hspace{-0.35mm} k \hspace{-0.35mm}<\hspace{-0.35mm} n$} (9);
%
        \draw [arrow, thick] (8) -- (10);
        \draw [arrow, thick] (9) -- (10);
%
        \draw [arrow, thick] (10) -- node[midway, fill = white, sloped] {\scriptsize $v \in \mu$} (11);
        \draw [arrow, thick] (10) -- node[midway, fill = white, sloped] {\scriptsize $v \not\in \mu$} (12);
%
        \draw [thick] (-1.46, -20.85) -- node[midway, fill = white, sloped] {\scriptsize $| \lambda | = 1$} (-8.5,-20.85);
        \draw [arrow, thick] (-8.5,-20.85) -- (6);
%
        \draw [thick] (15.5,-6.35) -- (15.5, -15.5);
        \draw [arrow, thick] (15.5, -15.5) -- (10);
    \end{tikzpicture}
    \end{align*}

\begin{proposition}\label{alwaysterm}
Strategy~\ref{thestrat} always terminates by correctly identifying the hidden ideal \(\mu\) in \(\lambda^0\), while using at most \(k^0\) positive queries.
\end{proposition}
\begin{proof}
If \(\lambda^0 = \varnothing\), then Strategy~\ref{thestrat} terminates by recognizing \(\mu = \varnothing\) and utilizes no queries. If \(|\lambda^0| = 1\), then Strategy~\ref{thestrat} terminates by recognizing either \(\mu = \varnothing\) or \(\mu = \lambda^0\), and utilizes \(1 \leq k^0\) queries. Now set \(|\lambda^0|=m>1\), and assume Strategy~\ref{thestrat} terminates by identifying any hidden ideal \(\mu^1\) in a poset \( \lambda^1\) utilizing at most \(k^1\) queries whenever \(|\lambda^1| < m\). 
We consider the application of Strategy~\ref{thestrat} to \(\lambda^0\). \\

\noindent{\em Case 1.} Assume \(k=1\). 
Then Strategy~\ref{thestrat} repeatedly queries nodes \(v\) of maximal height until either (a) a positive query results, in which case  Strategy~\ref{thestrat} correctly indicates that \(\mu = \lambda^0_{\preccurlyeq v}\), or; (b) Strategy~\ref{thestrat} terminates after (negatively) querying every node in \(\lambda^0\), and correctly indicates that \(\mu = \varnothing\).\\

\noindent{\em Case 2.} Assume \(k >1\). Then Strategy~\ref{thestrat} queries some node \(v \in \lambda\) of height greater than 1. We have subcases:
\begin{enumerate}[(i)]
\item Assume \(v \in \mu\). Set \(\lambda^1:=\lambda_{\succcurlyeq v}\), \(\mu^1:=\mu \cap \lambda_{\succcurlyeq v}\) and \(k^1:=k-1\). Then Strategy~\ref{thestrat} continues by running the algorithm for \((\lambda^1, k^1)\), and by induction assumption, terminates by correctly identifying \(\mu^1\) as some ideal \(\lambda^1_{\succcurlyeq w}\) in \(\lambda^1\) while using at most \(k^1\) queries. But this implies that \(\mu = \lambda^0_{\succcurlyeq w}\), and thus Strategy~\ref{thestrat} correctly identifies \(\mu\) using at most \(1+k^1=k\) positive queries.

\item Assume \(v \notin \mu\). Set \(\lambda^1:=\lambda_{\not\succcurlyeq v}\), \(\mu^1:=\mu \cap \lambda_{\not\succcurlyeq v} = \mu\), and \(k^1 :=k\). Then Strategy~\ref{thestrat} continues by running the algorithm for \((\lambda^1, k^1)\), and by induction assumption, terminates by correctly identifying \(\mu^1\) as some ideal \(\lambda^1_{\succcurlyeq w}\) (or \(\varnothing\)) in \(\lambda^1\) while using at most \(k\) queries. But this implies that \(\mu = \lambda^0_{\succcurlyeq w}\) (or \(\varnothing\)), and thus Strategy~\ref{thestrat} correctly identifies \(\mu\) using at most \(k\) positive queries.
\end{enumerate}
Thus in any case, Strategy~\ref{thestrat} terminates by correctly identifying \(\mu\) while using at most \(k\) positive queries.
\end{proof}

\subsection{Example}\label{exsec}
\begin{figure}[h]
\begin{align*}
        \begin{tikzpicture}[scale = 0.5]
            \draw(8,1)--(4,3);
            \draw(8,1)--(7,3);
            \draw(7,3)--(7,9);
            \draw(8,1)--(9,3) [ultra thick, bittersweet];
            \draw(8,1)--(12,3) [ultra thick, bittersweet];
            \draw(4,3)--(3,5);
            \draw(4,3)--(5,5);
            \draw(7,3)--(6,5);
            \draw(7,3)--(8,5);
            \draw(9,3)--(10,5);
            \draw(12,3)--(12,5) [ultra thick, bittersweet];
            \draw(12,3)--(14,5);
            \draw(14,5)--(14,7);
            \draw(14,5)--(15,7);
            \draw(3,5)--(2,7);
            \draw(3,5)--(3,7);
            \draw(3,5)--(4,7);
            \draw(5,5)--(6,7);
            \draw(9,3)--(8,5);
            \draw(2,7)--(1,9);
            \draw(2,7)--(2,9);
            \draw(2,7)--(3,9);
            \draw(6,7)--(5,9);
            \draw(6,7)--(6,9);
            \draw(6,7)--(7,9);
            \draw(3,9)--(9,11);
            \draw(12,9)--(9,11);
            \draw(6,9)--(9,11);
            \draw(10,9)--(9,11);
            \draw(12,5)--(9,7) [ultra thick, bittersweet];
            \draw(9,7)--(7,9);
            \draw(9,7)--(9,9);
            \draw(9,7)--(10,9) [ultra thick, bittersweet];
            \draw(9,3)--(9,7) [ultra thick, bittersweet];
            \draw(12,5)--(14,9);
            \draw(15,7)--(14,9);
            \draw(15,7)--(15,9);
            \draw(15,7)--(16,9);
            \draw(12,5)--(12,7);
            \draw(12,7)--(12,9);
            \draw(12,7)--(13,9);
            \draw(4,7)--(3,9);
%
            \draw[thick, bittersweet, fill = bittersweet] (8,1) circle(1.5mm);
            \draw[thick, black, fill = black] (4,3) circle(1.5mm);
            \draw[thick, black, fill = black] (7,3) circle(1.5mm);
            \draw[thick, black, fill = black] (15,9) circle(1.5mm);
            \draw[thick, black, fill = black] (12,7) circle(1.5mm);
            \draw[thick, black, fill = black] (12,9) circle(1.5mm);
            \draw[thick, black, fill = black] (13,9) circle(1.5mm);
            \draw[thick, black, fill = black] (15,9) circle(1.5mm);
            \draw[thick, black, fill = black] (14,7) circle(1.5mm);
            \draw[thick, black, fill = black] (14,5) circle(1.5mm);
            \draw[thick, black, fill = black] (15,7) circle(1.5mm);
            \draw[thick, black, fill = black] (14,9) circle(1.5mm);
            \draw[thick, bittersweet, fill = bittersweet] (12,3) circle(1.5mm);
            \draw[thick, bittersweet, fill = bittersweet] (12,5) circle(1.5mm);
            \draw[thick, bittersweet, fill = bittersweet] (9,3) circle(1.5mm);
            \draw[thick, black, fill = black] (10,5) circle(1.5mm);
            \draw[thick, black, fill = black] (6,5) circle(1.5mm);
            \draw[thick, black, fill = black] (8,5) circle(1.5mm);
            \draw[thick, black, fill = black] (5,5) circle(1.5mm);
            \draw[thick, black, fill = black] (3,5) circle(1.5mm);
            \draw[thick, black, fill = black] (2,7) circle(1.5mm);
            \draw[thick, black, fill = black] (1,9) circle(1.5mm);
            \draw[thick, black, fill = black] (2,9) circle(1.5mm);
            \draw[thick, black, fill = black] (3,9) circle(1.5mm);
            \draw[thick, black, fill = black] (3,7) circle(1.5mm);
            \draw[thick, black, fill = black] (4,7) circle(1.5mm);
            \draw[thick, black, fill = black] (9,11) circle(1.5mm);
            \draw[thick, black, fill = black] (5,9) circle(1.5mm);
            \draw[thick, black, fill = black] (6,7) circle(1.5mm);
            \draw[thick, black, fill = black] (6,9) circle(1.5mm);
            \draw[thick, black, fill = black] (16,9) circle(1.5mm);
            \draw[thick, black, fill = black] (7,9) circle(1.5mm);
            \draw[thick, bittersweet, fill = bittersweet] (9,7) circle(1.5mm);
            \draw[thick, black, fill = black] (9,9) circle(1.5mm);
            \draw[thick, bittersweet, fill = bittersweet] (10,9) circle(1.5mm);
%
            \draw[thick, cyan] (5,5) circle(3mm);
    \node[label=left:\textcolor{cyan}{\scriptsize $v_1$}] at (5,5) {};
            \draw[thick, cyan] (4,3) circle(3mm);
    \node[label=left:\textcolor{cyan}{\scriptsize $v_2$}] at (4,3) {};
            \draw[thick, cyan] (12,3) circle(3mm);
    \node[label=right:\textcolor{cyan}{\scriptsize $v_3$}] at (12,3) {};
            \draw[thick, cyan] (9,7) circle(3mm);
    \node[label=right:\textcolor{cyan}{\scriptsize $v_4$}] at (9,7) {};
            \draw[thick, cyan] (10,9) circle(3mm);
    \node[label=right:\textcolor{cyan}{\scriptsize $v_6$}] at (10,9) {};
            \draw[thick, cyan] (9,9) circle(3mm);
    \node[label=left:\textcolor{cyan}{\scriptsize $v_5$}] at (9,9) {};
        \end{tikzpicture}
    \end{align*}
    \caption{The poset \(\lambda^0\) of degree 4 and height 6, with hidden ideal \(\mu\) highlighted in orange, as described in Example~\ref{exsec}. The nodes \(v_1, \ldots, v_6\) queried during application of Strategy~\ref{thestrat} are circled in cyan.}
            \label{BigEx}
            \end{figure}
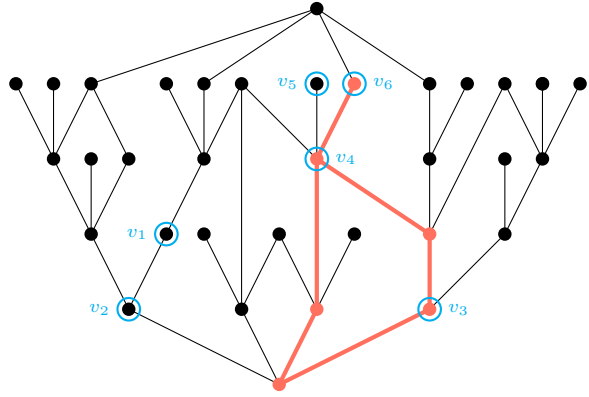

In Figure~\ref{BigEx} we display a poset \(\lambda^0\) of degree 4 and height 6. A hidden ideal \(\mu\) is shown highlighted in orange. We now show a step-by-step application of Strategy~\ref{thestrat} which identifies \(\mu\) while using at most \(k^0 = 3\) positive queries.

    \noindent\textit{Query 1:} We have \(\textup{ht}(\lambda) = \textup{ht}(\lambda^0) = 6\) and \(k=k^0 = 3\), so following Strategy~\ref{thestrat}, we query the node \(v_1\) of height \(\lceil \frac{6+1}{3}\rceil = 3\) that resides in a maximal chain in \(\lambda\). Since \(v_1 \notin \mu\), we set \(\lambda:= \lambda_{\not \succcurlyeq v_1}\), effectively erasing all nodes that dominate \(v_1\). \\

    \noindent\textit{Query 2:} We now have \(\textup{ht}(\lambda) = 5\) and \(k=3\), so following Strategy~\ref{thestrat}, we query the node \(v_2\) of height \(\lceil \frac{5+1}{3}\rceil = 2\) that resides in a maximal chain in \(\lambda\). Since \(v_2 \notin \mu\), we set \(\lambda:= \lambda_{\not \succcurlyeq v_2}\), effectively erasing all nodes that dominate \(v_2\). \\

    \noindent\textit{Query 3:} We still have \(\textup{ht}(\lambda) = 5\) and \(k=3\), so following Strategy~\ref{thestrat}, we query the node \(v_3\) of height \(\lceil \frac{5+1}{6}\rceil = 2\) that resides in a maximal chain in \(\lambda\). Since \(v_3 \in \mu\), we set \(\lambda:= \lambda_{ \succcurlyeq v_3}\), effectively erasing all nodes that do not dominate \(v_3\). Having used a positive query, we set \(k:=3-1 = 2\). \\

\noindent\textit{Query 4:} We now have \(\textup{ht}(\lambda) = 4\) and \(k=2\), so following Strategy~\ref{thestrat}, we query the node \(v_4\) of height \(\lceil \frac{4+1}{2}\rceil = 3\) that resides in a maximal chain in \(\lambda\). Since \(v_4 \in \mu\), we set \(\lambda:= \lambda_{ \succcurlyeq v_4}\), effectively erasing all nodes that do not dominate \(v_4\). Having used a positive query, we set \(k:=2-1 = 1\). \\

\noindent\textit{Query 5:} We now have \(\textup{ht}(\lambda) = 2\) and \(k=1\), so following Strategy~\ref{thestrat}, we query the node \(v_5\) of height \(2\). Since \(v_5 \notin \mu\), we set \(\lambda:= \lambda_{ \not\succcurlyeq v_5}\), effectively erasing \(v_5\).\\

\noindent\textit{Query 6:} We now have \(\textup{ht}(\lambda) = 2\) and \(k=1\), so following Strategy~\ref{thestrat}, we query the node \(v_6\) of height \(2\). Since \(v_6 \in \mu\), we have \(\mu= \lambda^0_{\preccurlyeq v_6}\), completing the identification strategy.

\section{Bounds for Strategy~\ref{thestrat}}\label{stratanalsec}
In this section we will establish bounds on the maximum number of total queries Strategy~\ref{thestrat} requires to identify a hidden ideal \(\mu\) in a poset \(\lambda\), while using at most \(k\) positive queries.
Recall from \S\ref{gensearsec} that we define a function
\begin{align*}
f_{k,\ell}(n) := \sum_{i=0}^{\lceil n/k \rceil-1} \ell^i + \sum_{j=1}^{k-1} \ell^{\lceil (n-j)/k\rceil}
\end{align*}
for all \(k \in \mathbb{N}\) and \(\ell,n \in \mathbb{Z}_{\geq 0}\). We take \(0^0 =1\) when \(\ell = 0\).
In some special cases, it is straightforward to check that this formula simplifies to:
\begin{align}\label{niceeq}
f_{k,\ell}(n) =
\begin{cases}
k & \textup{if } \ell = 0, n=1;\\
\lceil \frac{n}{k} \rceil + k -1 & \textup{if } \ell = 1;\\
n\ell - \ell + k - n + 1 & \textup{if }n \leq k.\\
\end{cases}
\end{align}

    \subsection{Bounds for Strategy~\ref{thestrat}} We now establish our first main result. The proof relies on some technical results on ceiling/floor functions and the behavior of the function \(f_{k,\ell}(n)\), which we relegate to the Appendix \S\ref{append}.
    \begin{theorem}\label{mainthmA}
        Let \(\ell, n \in \mathbb{Z}_{\geq 0}\), and let $\lambda \in \mathcal{P}(\ell,n)$. Let \(\mu \subseteq \lambda\) be a hidden ideal in \(\lambda\). For \(k \in \mathbb{N}\), Strategy~\ref{thestrat} identifies the hidden ideal $\mu$ using at most \(k\) positive queries and at most \(f_{k,\ell}(n)\) total queries.
    \end{theorem}

    \begin{proof}
    By Proposition~\ref{alwaysterm}, we have that Strategy~\ref{thestrat} identifies \(\mu\) using at most \(k\) positive queries, so we now work on establishing the bound \(f_{k,\ell}(n)\) on total queries.
        We go by induction on $|\lambda|$.\\

        \noindent\textit{Base Case:} Assume $|\lambda| = 1$. Then $\ell =0, n = 1$. Then Strategy~\ref{thestrat} identifies \(\mu\) with one query, and we have \(f(0,1,k) = k \geq 1\), as required.\\

        \noindent\textit{Induction Step:} Assume \(\lambda \in \mathcal{P}(\ell,n)\), $|\lambda| > 1$, and the claim holds for all $\lambda' \in \mathcal{P}(\ell',n')$ with $|\lambda'| < |\lambda|$. Say $\lambda \in \mathcal{P}(\ell,n)$. Note that since \(|\lambda|>1\), we have \(\ell \geq 1\) and \(n\geq 2\). Following Strategy~\ref{thestrat}, we have three cases to consider.\\
        
\noindent{\em Case 1.}
            Assume $k=1$. Following Strategy~\ref{thestrat}, we repeatedly query nodes of $\lambda$ in a top-down fashion until \(\mu\) is identified upon the first positive query. So at worst, Strategy~\ref{thestrat} requires $|\lambda|$ queries. But by Lemma~\ref{new_lem}, we have 
            \begin{align*}
|\lambda| \leq \sum_{i=0}^{n-1} \ell^i = \sum_{i=0}^{\lceil n/1 \rceil-1} \ell^i = f_{1,\ell}(n),
            \end{align*}
            as required.\\

\noindent{\em Case 2.}
Assume \(k \geq n > 1\). 
Following Strategy~\ref{thestrat}, we set \(v \in \lambda\) to be a node of height 2 in a maximal chain in \(\lambda\), and query \(v\). We have two possible outcomes:\\

\noindent{\em Case 2.1.} Assume \(v \notin \mu\). Then \(|\lambda_{\not \succcurlyeq v}| = 1\), and following Strategy~\ref{thestrat}, we identify \(\mu\) after one additional query of the base node in \(\lambda\), for a grand total of two queries. But then by (\ref{niceeq}) we have
\begin{align*}
f_{k,\ell}(n) = (n-1)\ell + k - n + 1 \geq \ell + 1 \geq 2,
\end{align*}
as required.\\

\noindent{\em Case 2.2.} Assume \(v \in \mu\). Then we have used up a positive query, and now have \(k-1\) positive queries remaining. Moreover, \(\lambda_{\succcurlyeq v} \in \mathcal{P}(\ell', n-1)\) for some \(\ell' \leq \ell\), and \(|\lambda_{\succcurlyeq v}| < |\lambda|\). Thus by induction assumption, Strategy~\ref{thestrat} identifies \(\mu\) in a total of less than or equal to \(f_{k-1,\ell}(n-1)+1\) queries. But by (\ref{niceeq}) we have
\begin{align*}
f_{k-1,\ell}(n-1)+1&= (n-2)\ell + k-n + 2 \leq (n-1)\ell + k-n + 1 = f_{k,\ell}(n),
\end{align*}
as required.\\

\noindent{\em Case 3.} Assume
            $n > k > 1$. Set $j = \left\lceil \frac{n+1}{k} \right\rceil$. Assume there are $t$ nodes of height $j$ in $\lambda: v_1,v_2,\dots, v_t$. Note that $t$ is at most $\ell^{j-1}$ by Lemma~\ref{new_lem}. We have three possible subcases: 
                \begin{enumerate}
                    \item[(1)] Strategy~\ref{thestrat} terminates with a solution in less than or equal to $t$ steps.
                    \item[(2)] Strategy~\ref{thestrat} queries some nodes $w_1, w_2, \dots, w_{r-1}$, and all return negative results (i.e. not in $\mu$) until $w_r$ returns a positive result, with $r\leq t$.
                    \item[(3)] Strategy~\ref{thestrat} queries $w_1, w_2, \dots, w_t$ and all return negative results.
                \end{enumerate}
            
            \noindent\textit{Case 3.1.} Assume Strategy~\ref{thestrat} terminates with a  solution in less than or equal to $t$ steps. Then we have
            \begin{align*}
                t \leq \ell^{j-1} = \ell^{\lceil (n+1)/k \rceil - 1} = \ell^{\lceil n-(k-1)/k \rceil} \leq \sum_{i=0}^{\lceil n/k \rceil -1} \ell^i + \sum_{t=1}^{k-1} \ell^{\lceil(n-t)/k \rceil} = f_{k,\ell}(n)
            \end{align*}
                as required.\\

            \noindent\textit{Case 3.2.} Assume Strategy~\ref{thestrat} queries nodes $w_1, w_2, \dots, w_{r-1}$, all with negative results (i.e. not in $\mu$) until $w_r$ returns with a positive result, with $r \leq t$. Set $n_0 = n$. Set $n_i = \textup{ht}(\lambda_{\not\succcurlyeq\{w_1, \dots, w_i\}})$ for $i = 1, \dots, r-1$. It follows from Strategy~\ref{thestrat} and the above assumptions that $\textup{ht}(w_i) = \left\lceil \frac{n_{i-1}+1}{k} \right\rceil$. Then, $\textup{ht}(w_r) = \left\lceil \frac{n_{r-1}+1}{k} \right\rceil$. Also, note $n = n_0 \geq n_1 \geq \cdots \geq n_{r-1}$. Since $w_r$ is a height $\left \lceil \frac{n_{r-1}+1}{k} \right\rceil$ node in a maximal chain in $\lambda_{\not\succcurlyeq \{w_1,\dots, w_{r-1}\}}$, we have
                \begin{align*}
                    \textup{ht}((\lambda_{\not\succcurlyeq\{w_1,\dots, w_{r-1}\}})_{\succcurlyeq w_r}) &= n_{r-1} - \left\lceil \frac{n_{r-1}+1}{k} \right\rceil + 1 
                    = \left\lfloor \frac{(n_{r-1}+1)(k-1)}{k} \right\rfloor\\
                    &\leq \left\lfloor \frac{(n+1)(k-1)}{k} \right\rfloor
                    = n - \left\lceil \frac{n+1}{k} \right\rceil + 1,
                \end{align*}
            where the second and third equalities follow from Lemma~\ref{RR123}.   
            Set $\lambda' = (\lambda_{\not\succcurlyeq\{w_1, \dots, w_{r-1}\}})_{\succcurlyeq w_r}$, and $n' = \textup{ht}\left( (\lambda_{\not\succcurlyeq\{w_1,\dots,w_{r-1}\}})_{\succcurlyeq w_r} \right)$. Since \(|\lambda'| < |\lambda|\), and \(\lambda' \in \mathcal{P}(\ell',n')\), we have by induction assumption that Strategy~\ref{thestrat} solves $\lambda'$ can in less than or equal to $f_{k-1,\ell}(n')$ queries. Thus Strategy~\ref{thestrat} solves $\lambda$ in less than or equal to \(r + f_{k-1,\ell}(n')\) queries. We have
             \begin{align*}
                    r + f_{k-1,\ell}(n')
                    &\leq 
                    \ell^{\left\lceil \frac{n+1}{k} \right\rceil - 1} + 
                    \sum_{i=0}^{\left\lceil \frac{n-\left\lceil \frac{n+1}{k} \right\rceil + 1}{k-1} \right\rceil-1} \ell^i
                    +
                    \sum_{t=1}^{k-2} \ell^{\left\lceil \frac{n-\left\lceil \frac{n+1}{k} \right\rceil + 1 - t}{k-1} \right\rceil}\\
                    &=
                    \ell^{\left\lceil \frac{n-(k-1)}{k} \right\rceil}
                    +
                    \sum_{i=0}^{\left\lceil \frac{\left\lceil \frac{n(k-1)}{k} \right\rceil }{k-1} \right\rceil-1} \ell^i
                    +
                    \sum_{t=1}^{k-2}\ell^{\left\lceil \frac{\left\lceil \frac{n(k-1)}{k} \right\rceil - t}{k-1} \right\rceil}
                    \\
                    &=
                    \ell^{\left\lceil \frac{n-(k-1)}{k} \right\rceil}
                    +
                    \sum_{i=0}^{\left\lceil \frac{n}{k}\right\rceil-1} \ell^i
                    +
                    \sum_{t=1}^{k-2} \ell^{\lceil \frac{n-t}{k} \rceil}\\
                    &= f_{k,\ell}(n).
            \end{align*}        
                where the first equality follows by Lemma \ref{lem_n-k+1}, and the second equality follows by Lemma \ref{lem_n_k_i}.\\
           \\
            \textit{Case 3.3.} Assume the algorithm queries $w_1, \dots, w_t$, and all return negative results. Then, note that $\textup{ht}(w_i) \leq j$ for all $i \in [1,t]$ and, by construction of the algorithm, each $w_i$ belongs to a maximal chain in $\lambda_{\not\succcurlyeq\{w_1, \dots, w_{i-1}\}}$. Thus, by Lemma \ref{lem_height<j}, we have $\textup{ht}(\lambda_{\not\succcurlyeq\{w_1, \dots, w_t\}}) < j$. So, after querying $w_1,\dots, w_t$, we continue by applying Strategy~\ref{thestrat} to $\lambda_{\not\succcurlyeq\{w_1,\dots,w_t\}}$. By induction assumption, Strategy~\ref{thestrat} solves $\lambda_{\not\succcurlyeq\{w_1,\dots,w_t\}}$ in less than or equal to $f_{k,\ell}(j-1)$  queries, and thus solves \(\lambda\) in less than or equal to \(t + f_{k,\ell}(j-1)\) queries.  By Lemma \ref{n/k=j}, we have
                \begin{align*}
                    t + f_{k,\ell}(j-1) \leq \ell^{j-1} + f_{k,\ell}(j-1) \leq f_{k,\ell}(n),
                \end{align*}
                which completes the induction step and the proof.
    \end{proof}

Now recall from \S\ref{HIprob} the value \(\boldsymbol{q}_k(\lambda)\) which records the {\em minimal} necessary queries to identify a hidden ideal \(\mu\) in a poset \(\lambda\) granted at most \(k\) positive queries. 

\begin{corollary}\label{corbound}
Let \(k \in \mathbb{N}\), and let \(\lambda\) be a pointed poset of degree \(\ell\) and height \(n\). Then we have
\begin{align*}
\boldsymbol{q}_k(\lambda) \leq f_{k,\ell}(n) = \sum_{i=0}^{\lceil n/k \rceil-1} \ell^i + \sum_{j=1}^{k-1} \ell^{\lceil (n-j)/k\rceil} \leq (k \ell) \ell^{n/k}.
\end{align*}
\end{corollary}
\begin{proof}
The first inequality follows from Theorem~\ref{mainthmA}, so we focus now on the second inequality.
If \(\ell = 1\), this is immediate from (\ref{niceeq}). If \(\ell \geq 2\) we have
\begin{align*}
f_{k,\ell}(n) &= \sum_{i=0}^{\lceil n/k \rceil-1} \ell^i + \sum_{j=1}^{k-1} \ell^{\lceil (n-j)/k\rceil}
= \frac{\ell^{\lceil n/k \rceil}-1}{\ell-1} + \sum_{j=1}^{k-1} \ell^{\lceil (n-j)/k\rceil}\\
&\leq \ell^{\lceil n/k \rceil} + \sum_{j=1}^{k-1} \ell^{\lceil n/k\rceil} = k \ell^{\lceil n/k \rceil} \leq k \ell^{ 1+ n/k} =( k\ell) \ell^{n/k},
\end{align*}
as required.
\end{proof}

\section{Asymptotic optimality of Strategy~\ref{thestrat}}\label{efficsec}
We first recall here some results from \cite{IJM} which will be useful. 
For any \(x \in \mathbb{Z}_{\geq 0}\), let
\(
T_k(x) := \sum_{i=1}^k {x \choose k},
\)
and define \(\tau_k(x)\) to be the smallest integer such that \(x \leq T_k(\tau_k(x))\). The following result from  \cite[Theorem~4.2]{IJM} gives lower bounds on the value of \(\boldsymbol{q}_k(\lambda)\).

\begin{theorem}\label{oldboundthm}
For all \(k \in \mathbb{N}\) and posets \(\lambda\), we have \(\tau_k(|\lambda|) \leq \boldsymbol{q}_k(\lambda) \).
\end{theorem}

Recall that \(\mathcal{T}_\ell(n)\) is the complete \(\ell\)-ary tree of height \(n\).

\begin{lemma}\label{mbound}
Fix \(k, \ell \in \mathbb{N}\). There exists \(m > 0\) such that for all \(n \in \mathbb{N}\), we have \(m \ell^{n/k} \leq \boldsymbol{q}_k(\mathcal{T}_\ell(n))\).
\end{lemma}
\begin{proof}
First, note that
\begin{align*}
T_k(x) = \sum_{i=1}^k {x \choose k} = \sum_{i=1}^k \frac{1}{k!} x(x-1)\cdots (x-k+1)
\end{align*}
is a degree \(k\) polynomial in \(x\), and that by definition \(T_k(x)\) is an increasing function on \(x\). Thus we may write
\(
T_k(x) = \sum_{i=0}^k a_i x^i
\)
for some coefficients \(a_0, \ldots, a_k \in \mathbb{R}\). Set \(u = |a_0| + \cdots + |a_k|\). Then for all \(x \in \mathbb{N}\) we have
\begin{align*}
T_k(x) = \sum_{i=0}^k a_i x^i \leq \sum_{i=0}^k |a_i| x^i \leq \sum_{i=0}^k |a_i| x^k = ux^k.
\end{align*}
Now define a function \(\Gamma_k(x) = \left( \frac{x}{A} \right)^{1/k}\), and set \(m= (\ell u)^{-1/k}\). Since \(T_k(x) \leq ux^k\) and \(T_k\) is increasing in \(x\), we have
\begin{align*}
T_k(\Gamma_k(x)) \leq \Gamma_k(ux^k) = \left( \frac{ux^k}{u}\right)^{1/k} = x \leq T_k(\tau_k(x)).
\end{align*}
which again implies that \(\Gamma_k(x) \leq \tau_k(x)\) since \(T_k\) is increasing. Then for  \(n \in \mathbb{N}\) we have
\begin{align*}
m\ell^{n/k} &= \frac{\ell^{n/k}}{(\ell u)^{1/k}} =
\frac{\ell^{(n-1)/k}}{u^{1/k}}=
\left( \frac{\ell^{n-1}}{u} \right)^{1/k}=\Gamma_k(\ell^{n-1})\\
&\leq \Gamma_k\left(\sum_{i=0}^{n-1} \ell^i\right) = \Gamma_k(|\mathcal{T}_\ell(n)|) \leq \tau_k(|\mathcal{T}_\ell(n)|) \leq \boldsymbol{q}_k(\mathcal{T}_\ell(n)),
\end{align*}
where we have used the fact that \(\Gamma_k\) is increasing for the first inequality, and Theorem~\ref{oldboundthm} for the last inequality. 
\end{proof}

Now we are finally in position to prove our second main result, which shows that Strategy~\ref{thestrat} performs asymptotically optimally on the family of \(\ell\)-ary trees as the height \(n\) grows.

\begin{theorem}\label{growthm}
Fix \(k,\ell \in \mathbb{N}\).
Then there exist \(m,M > 0\) such that
\begin{align*}
m \ell^{\,n/k} \leq \boldsymbol{q}_{k}(\mathcal{T}_\ell(n)) \leq f_{k,\ell}(n) \leq M \ell^{\,n/k},
\end{align*}
for all \(n \in \mathbb{N}\). 
Thus \(\boldsymbol{q}_{k}(\mathcal{T}_\ell(n)) = \Theta(\ell^{\,n/k}) = f_{k,\ell}(n)\) as functions of \(n\).
\end{theorem}
\begin{proof}
By Lemma~\ref{mbound} there exists \(m>0\) such that \(m \ell^{n/k} \leq \boldsymbol{q}_k(\mathcal{T}_\ell(n))\). Taking \(M = k\ell\), we have by Corollary~\ref{corbound} that \(\boldsymbol{q}_k(\mathcal{T}_\ell(n)) \leq f_{k,\ell}(n) \leq M \ell^{n/k}\) since \(\mathcal{T}_\ell(n) \in \mathcal{P}(\ell,n)\), completing the proof.
\end{proof}

    \section{Appendix}\label{append}
    In this section, we will establish the results on the function \(f_{k,\ell}(n)\) and ceiling/floor functions that are used in our proof of Theorem~\ref{mainthmA}.
\begin{lemma}\label{RR123}
Let \(m,k \in \mathbb{N}\). Then we have \(m- \lceil \frac{m+1}{k}\rceil + 1 = \lfloor \frac{(m+1)(k-1)}{k}\rfloor\).
\end{lemma}
\begin{proof}
We have
 \begin{align*}
                   m - \left\lceil \frac{m+1}{k} \right\rceil + 1 
                    &= -\left(-m + \left\lceil \frac{m + 1}{k} \right\rceil -1 \right)
                    = - \left\lceil -m + \frac{m + 1}{k} -1 \right\rceil\\
                    &= - \left\lceil \frac{(m + 1)(1-k)}{k} \right\rceil
                    = \left\lfloor \frac{(m+1)(k-1)}{k} \right\rfloor,
                \end{align*}
                as required.
\end{proof}

    \begin{lemma} \label{lem_n-k+1}
        Let $n,k \in \mathbb{N}$. Then we have
            \begin{align*}
                n - \left\lceil \frac{(k-1)n}{k} \right\rceil = \left\lceil \frac{n-k+1}{k} \right\rceil = \left\lceil \frac{n+1}{k} \right\rceil -1. 
            \end{align*}
        
    \end{lemma}

    \begin{proof}
        The second equality is obvious. We prove the first equality. By the division algorithm we have $n = kq +r$ for some $r \in [0, k -1]$ and $q \in \mathbb{Z}$. Then we have
        \begin{align*}
            n- \left\lceil\frac{(k-1)n}{k} \right\rceil &=
            (kq+r) - \left\lceil \frac{(k-1)(kq+r)}{k} \right \rceil = (kq+r) - \left\lceil \frac{k^2q+kr-kq-r}{k} \right\rceil\\
            &= kq+r - \left\lceil kq+r-q-\frac{r}{k} \right\rceil
            = kq+r - kq - r + q + \left\lceil \frac{-r}{k} \right\rceil\\
            &= q + \left\lceil \frac{-r}{k} \right\rceil = q.
        \end{align*}
        On the other hand, we have
        \begin{align*}
        \left \lceil \frac{n-k+1}{k}\right \rceil &=
            \left\lceil \frac{(kq+r) - (k-1)}{k} \right\rceil = \left\lceil \frac{kq-k+r+1}{k} \right\rceil\\
            &= \left\lceil (q-1) + \frac{r+1}{k} \right\rceil
            = q - 1 + \left\lceil \frac{r+1}{k} \right\rceil
            = q-1 + 1 = q.
        \end{align*}
Therefore, $n - \left\lceil \frac{(k-1)n}{k} \right\rceil = \left\lceil \frac{n-(k-1)}{k} \right\rceil$, as required.
    \end{proof}

    \begin{lemma} \label{lem_n_k_i}
    Let $n,k \in \mathbb{N}$, \(k \geq 2\), and let $i \in [0, k - 2]$. Then we have
        \begin{align*}
            \left\lceil \frac{\left\lceil \frac{nk-n}{k} \right\rceil - i}{k - 1} \right\rceil = \left\lceil \frac{n - i}{k} \right\rceil.
        \end{align*}
    \end{lemma}
    
    \begin{proof}
        Write $n = kq+r$ for some $r \in [0, k -1]$ and $q \in \mathbb{Z}$. Then we have
            \begin{align*}
            \left\lceil \frac{\left\lceil \frac{nk-n}{k} \right\rceil - i}{k - 1} \right\rceil &=
                \left\lceil \frac{\left\lceil \frac{(k-1)(kq+r)}{k} \right\rceil - i}{k-1} \right\rceil = \left\lceil \frac{\left\lceil \frac{k^2q + kr -kq - r}{k} \right\rceil - i}{k-1} \right\rceil\\
                &= \left\lceil \frac{(kq+r-q+\left\lceil \frac{-r}{k} \right\rceil) - i}{k-1} \right\rceil
                =\left\lceil \frac{(kq+r-q+0) - i}{k-1} \right\rceil\\
                &=
                \left\lceil q + \frac{r-i}{k-1} \right\rceil = q + \left\lceil \frac{r-i}{k-1} \right\rceil\\
                &= \begin{cases}
                q & \textup{if } r \leq i\\
                q+1 & \textup{if }r > i.
                \end{cases}
            \end{align*}
            On the other hand, we have
            \begin{align*}
            \left \lceil \frac{n-i}{k} \right \rceil &= 
                \left\lceil \frac{kq+r-i}{k} \right\rceil 
                = \left\lceil q + \frac{r-i}{k} \right\rceil
                = q + \left\lceil \frac{r-i}{k} \right\rceil
                = 
                \begin{cases}
                q  & \textup{if }r \leq i \\
                q+1 & \textup{if }r>i.
                \end{cases}
            \end{align*}
            Therefore, $\left\lceil \frac{\left\lceil \frac{(k-1)n}{k} \right\rceil - i}{k - 1} \right\rceil = \left\lceil \frac{n - i}{k} \right\rceil$, as required.
    \end{proof}

    \begin{lemma} \label{n/k=j}
        Assume \(n > k \geq 2\), \(\ell \geq 1\), and set $j = \left\lceil \frac{n+1}{k} \right\rceil$. Then we have:
            \begin{align*}
                \ell^{j-1} + f_{k,\ell}(j-1) \leq f_{k,\ell}(n).
            \end{align*}
    \end{lemma}

    \begin{proof}
       
We first consider the case \(\ell = 1\). 
Note that \(n + 1 \geq k+2\), and so 
\begin{align*}
k^2 \leq k^2 + k -2 = (k+2)(k-1) \leq (n+1)(k-1) = nk - n + k -1.
\end{align*}
Therefore we have \(k^2   - nk  +n - k+1 \leq 0\), and so 
\begin{align*}
0 \geq \left\lceil \frac{k^2   - nk  +n - k+1}{k} \right \rceil =  \left\lceil \frac{n +1}{k} \right \rceil -1 + k -n = j-1 + k - n.
\end{align*}
Thus \(k+j-1 \leq n\), and so 
\begin{align*}
1 + \left \lceil \frac{j-1}{k} \right \rceil =  \left \lceil \frac{k+ j-1}{k} \right \rceil \leq  \left \lceil \frac{n}{k} \right \rceil,
\end{align*}
which implies in view of (\ref{niceeq}) that
\begin{align*}
1^{j-1} + f_{k,1}(j-1) = 1 + \left \lceil \frac{j-1}{k} \right \rceil + k-1 \leq  \left \lceil \frac{n}{k} \right \rceil + k-1 = f_{k,1}(n),
\end{align*}
completing the proof in the case \(\ell = 1\). From now on we assume \(\ell \geq 2\).\\

\noindent{\em Case 1.} Assume \(\lceil \frac{n}{k} \rceil = \lceil \frac{n+1}{k} \rceil = j\). Note that since  that \((k-1)(j-1) \geq 0\), we have \(-1 \geq \frac{j-jk-1}{k}\), and thus 
\(
0 > \lceil  \frac{j-jk-1}{k} \rceil = \lceil  \frac{j-1}{k} \rceil  - j
\). Therefore \(j-2 \geq \lceil \frac{j-1}{k} \rceil - 1\). Note that since \(n \geq \lceil \frac{n}{k} \rceil = j > j-1\), we have \(\lceil \frac{n-u}{k} \rceil \geq \lceil \frac{j-1-u}{k} \rceil\) for all \(u \in [1,k-1]\). Therefore
\begin{align*}
f_{k,\ell}(n) &=
\ell^{j-1} + \sum_{t = 0}^{j-2} \ell^t + \sum_{u = 1}^{k-1} \ell^{\lceil (n-u)/k\rceil}
\geq \ell^{j-1} +  \sum_{t = 0}^{\lceil (j-1)/k\rceil-1} \ell^t + \sum_{u=1}^{k-1} \ell^{\lceil (j-1-u)/k \rceil}
= \ell^{j-1} + f_{k,\ell}(j-1),
\end{align*}
as required.\\

\noindent{\em Case 2.} Assume \(k=2\), and \(\lceil \frac{n}{2} \rceil < \lceil \frac{n+1}{2}\rceil = \lceil \frac{n}{2} \rceil + 1\). It follows then that \(n\equiv 0 \pmod{4}\) or \(n \equiv 2 \pmod 4\). We approach these subcases separately.\\

\noindent{\em Case 2.1.} Assume \(n \equiv 0 \pmod 4\). Then \(n = 4m\) for some \(m \geq 1\). Then \(j = \lceil \frac{4m+1}{2} \rceil = 2m+1\). Then we have
\begin{align*}
f(\ell, n,k) &= \sum_{t=0}^{2m-1} \ell^t + \ell^{2m} \geq \ell^{2m} + \sum_{t=0}^{m} \ell^t = \ell^{j-1} + f(\ell,j-1,k),
\end{align*}
as required.\\

\noindent{\em Case 2.2.} Assume \(n \equiv 2 \pmod 4\). Then \(n = 4m+2\) for some \(m \geq 1\), and \(j = \lceil \frac{4m+3}{2}\rceil = 2m+2\). Then we have
\begin{align*}
f(\ell, n,k) = \sum_{t=0}^{2m} \ell^t + \ell^{2m+1} \geq \ell^{2m+1 }+ \sum_{t=0}^{m} \ell^t + \ell^m = \ell^{j-1} + f(\ell, j-1,k),
\end{align*}
as required.\\

\noindent{\em Case 3.} Finally, assume \(k \geq 3\) and   \(\lceil \frac{n}{k} \rceil < \lceil \frac{n+1}{k} \rceil = \lceil \frac{n}{k} \rceil + 1\). Then we have \(n=km\) and  \(j =m+ 1\) for some \(m\geq 2\). Since \((k-1)m \geq k-1\), we have \(m \geq \frac{m+k-1}{k}\). Therefore we have
\begin{align*}
\ell^{\lceil (n-1)/k \rceil} &=\ell^{\lceil m - 1/k \rceil} = \ell^m = \ell^{\lceil m \rceil} \geq \ell^{\lceil (m+k-1)/k\rceil} = \ell^{\lceil (m-1)/k \rceil + 1} = \ell \cdot \ell^{\lceil (m-1)/k \rceil } \notag \\
&\geq   2  \ell^{\lceil (m-1)/k \rceil }
=\ell^{\lceil (m-1)/k \rceil } + \ell^{\lceil (m-1)/k \rceil } \geq \ell^{\lceil (m-1)/k \rceil } + \ell^{\lceil (m-2)/k \rceil }.
\end{align*}
and
\begin{align*}
\ell^{\lceil (n-2)/k \rceil} &=\ell^{\lceil m - 2/k \rceil} = \ell^m = \ell^{j-1}.
\end{align*}
Thus we have
\begin{align*}
f_{k,\ell}(n) &= \sum_{i = 0}^{m-1} \ell^i + \ell^{\lceil (n-1)/k\rceil} +  \ell^{\lceil (n-2)/k\rceil} +  \sum_{t=3}^{k-1}\ell^{\lceil (n-t)/k \rceil}  \\
&\geq   \sum_{i = 0}^{\lceil \frac{m}{k}\rceil-1} \ell^i+ \left(\ell^{\lceil (m-1)/k \rceil} +  \ell^{\lceil (m-2)/k \rceil}\right) +  \ell^{j-1} + \sum_{t=3}^{k-1} \ell^{\lceil (m-t)/k \rceil}\\
&= \ell^{j-1} +  f_{k,\ell}(j-1),
\end{align*}
as required, completing the proof.
\end{proof}

\end{document}